\documentclass[a4paper]{article}
\usepackage{amssymb}
\newtheorem{Th}{Theorem}

\newtheorem{Lma}{Lemma}[section]

\newcommand{\be}{\begin{equation}}
\newcommand{\ee}{\end{equation}}
\newcommand{\R}{\mathbb{R}}
\newcommand{\N}{\mathbb{N}}
\newcommand{\C}{\mathbb{C}}
\newcommand{\Z}{\mathbb{Z}}

\newcommand{\reset}{\setcounter{equation}{0}\setcounter{Th}{0}\setcounter{Prop}{0}\setcounter{Co}{0}
\setcounter{Lm}{0}\setcounter{Rm}{0}}

\def\ti{\tilde}
\def\lf{\left}
\def\rg{\right}

\def\al{\alpha}
\def\la{\lambda}

\def\ep{\varepsilon}
\def\ds{\displaystyle}
\def\ov{\overline}
\def\Om{\Omega}
\def\om{\omega}
\def\p{\partial}

\begin{document}
\title{Conservation laws for conformal invariant variational problems.}
\author{Tristan Rivi\`ere\footnote{Department of Mathematics, ETH Zentrum,
CH-8093 Z\"urich, Switzerland.}}
\maketitle
{\bf Abstract :} {\it We succeed in writing 2-dimensional conformally invariant
non-linear elliptic PDE (harmonic map equation, prescribed mean curvature
equations...etc) in divergence form.  This divergence free quantities generalize
to target manifolds without symmetries the well known conservation laws for harmonic maps into homogeneous spaces.
From this form we can recover, without the use
of moving frame,  all the classical regularity results known for 2-dimensional conformally invariant
non-linear elliptic PDE (see \cite{Hel}) . It enable us also to establish new results.
In particular we solve a conjecture by E.Heinz asserting that the solutions to the precribed  bounded mean
curvature equation in arbitrary manifolds are continuous. }

\section{Introduction}
\reset
The abscence of possible applications of the maximum principle to solutions to non-linear
elliptic systems reduces drastically the tools available for answering questions regarding
the symmetry, the uniqueness or the regularity of these solutions. In such an
impoverishment of the available technics while passing from 
scalar PDE to systems, the search for conservation laws is, however, one of the remaining relevant strategy to adress these questions.
Harmonic maps into spheres give a good illustration of the efficiency of conservation laws in this setting.
An harmonic map $u$ from the n-dimensional unit ball $B^n$ into the unit sphere $S^{m-1}$ of ${\R}^m$
is a $W^{1,2}(B^n,{\R}^m)$ map (maps in $L^2$ whose first derivatives are also in $L^2$)
which takes value almost everywhere in the sphere $S^{m-1}$ and which solves the following PDE
\be
\label{I.1}
\begin{array}{l}
-\Delta u=u\,|\nabla u|^2
\end{array}
\ee
where $\Delta$ is the negative laplacian in ${\R}^n$ : $\Delta=\sum_{i=1}^n\frac{\p^2}{\p x_i^2}$. They are the critical points
of the Dirichlet energy $E(u)=\int_{B^n}|\nabla u|^2\, dx_1\cdots dx_n$ for all perturbations of the form $u_t=u-t\phi/|u-t\phi|$ where $\phi$
is an arbitrary compactly supported smooth map from $B^n$ into ${\R}^m$.  Because of the conformal invariance of the
Dirichlet energy $E$ in 2 dimension ($E(u\circ\varphi)=E(u)$ for arbitrary $u$ in $W^{1,2}({\R}^2,{\R}^m)$ and arbitrary conformal map $\varphi$ from ${\R}^2$ into ${\R}^2$), the harmonic map equation (\ref{I.1}) is conformal invariant 
in 2 dimension : if $u$ is a solution to (\ref{I.1}) in $W^{1,2}(B^2,{\R}^m)$, the composition with an arbitrary conformal map
$\varphi$ : $u\circ \varphi$ is again a solution to (\ref{I.1}). The conformal dimension 2 is also the critical dimension 
for (\ref{I.1}) : The left-hand-side of (\ref{I.1}) for a $W^{1,2}$ solution is in $L^1$, therefore a solution has a laplacian in $L^1$
which is the borderline case in 2 dimension which "almost" ensure that the first derivatives  are in $L^2$  (using standard estimates on Riesz potential \cite{Ste}).  So, in some sense, by inserting the $W^{1,2}$ bound assumption in the non-linearity we are almost back on our feet by bootstraping this regularity information in the linear part of the equation. None of the two sides,
linear and non-linear, of the equation is really dominant : the equation is critical.
Among the fundamental analysis issues regarding equation (\ref{I.1}) are {\bf 1)} the regularity of solution 
in conformal 2-dimension and {\bf 2)} the passage to the limit in the equation for sequences of solutions
having bounded $E$ energy. Both questions were solved by the introduction of the following conservation laws
discovered by J. Shatah  (\cite{Sha}) : $u$ is a solution to (\ref{I.1}) in $W^{1,2}$ if and only if the following holds
\be
\label{I.2}
\begin{array}{l}
\forall i,j\in\{1\cdots m\}\quad\quad div\lf(u^i\nabla u^j-u^j\nabla u^i\rg)=0\quad.
\end{array}
\ee
The cancellation of these divergences can be interpreted by the mean of Noether theorem using the symmetries
of the target $S^{m-1}$ (see Helein's book \cite{Hel}).  Using this form it becomes straightforward to aswer to the question 
{\bf 2)} using the compactness of the embedding of $W^{1,2}$ into $L^2$ (Rellich Kondrachov embedding Theorem).
The answer to question {\bf 1)} : the fact that $W^{1,2}$ solutions to (\ref{I.1}) are real analytic was established by F.H\'elein
(see \cite{Hel})  starting  again from the conservation laws (\ref{I.2}). The main step was to prove the continuity of the solution
since, by classical results in \cite{HiW}, \cite{LaU} and \cite{Mo}, continuous solutions are real analytic. Using the conservation laws (\ref{I.2}) and the fact that $\sum_j u^j\nabla u^j=0$, F. H\'elein wrote the equations (\ref{I.1}) in the following way :
\be
\label{I.3}
\begin{array}{l}
\ds-\Delta u^i=\sum_{j=1} u^i\, \nabla u_j\cdot\nabla u^j=\sum_{j=1}^m\lf[u^i\, \nabla u_j-u_j\,\nabla u^i\rg]\cdot\nabla u^j\\[5mm]
\ds \quad=\sum_{j=1}^n curl B^{i}_j\cdot\nabla u^j
\end{array}
\ee
The existence of $B^{i}_j$ in $W^{1,2}_{loc}$ solving $curl B^{i}_j=u^i\, \nabla u_j-u_j\,\nabla u^i$ is given by the classical theory of elliptic operators. The product curl-grad in the right-hand-side of (\ref{I.3}) has in fact some additional regularity
than being simply in $L^1$ and the inverse by the laplace operator of such a product is continuous. This special
phenomenon that we recall in the appendix was first observed in a particular case in \cite{We} by H.Wente and was proved
in it's full generality by H.Brezis and J.M. Coron in \cite{BrC} extending  Wente's argument and independantly, using a quite different approach, by L.Tartar in \cite{Ta1}. Later on, the product curl-grad was observed to be in the local  Hardy  space
${\mathcal H}^1_{loc}$, smaller than $L^1$,  by  R. Coifman, P.L. Lions, Y.Meyer and S.Semmes in \cite{CLMS} following the work of S. M\"uller
\cite{Mu} where this result was obtained under some sign assumption on the product. Among the special features of distributions
in ${\mathcal H}^1_{loc}$ is the "nice" behavior of this space with respect to Calderon-Zygmund operators, in particular
the inverse of such a distribution by the Laplace operator are in $W^{2,1}$ which embedds in $C^0$ in 2 dimension.
Observing that the non-linearity of the harmonic map equation is in ${\mathcal H}^1_{loc}$ gives not only another approach
to conclude that solutions to (\ref{I.1}) are smooth in 2 dimension but also permits to establish the estimate
\be
\label{I.4}
\int_{O}|\nabla^2 u|\, dx_1\cdots dx_n<+\infty\quad,
\ee
for any solution to (\ref{I.1}),  for $n$ arbitrary and where $O$ is an arbitrary open subset with closure in $B^n$. This estimate
happens to play a crucial role for establishing energy quantization results as described in \cite{LiR}.

It is now natural to try to understand to which extend the above results are still valid when we are 
considering $W^{1,2}$ harmonic maps taking values in an arbitrary submanifold of ${\R}^m$. What about questions
{\bf 1)}, {\bf 2}) or  what about the validity of the estimate (\ref{I.4}) in this general setting ?
Let then, $N^k$ be a $C^2$ k-dimensional submanifold of ${\R}^m$. Denote $\pi_{N}$ the $C^1$ orthogonal projection on $N^k$ defined in a small tubular neighborhood of $N^k$ in ${\R}^m$ which assigns
to each point in this neighborhood the nearest point on $N$. We denote by $W^{1,2}(B^n,N^k)$
the subset of $W^{1,2}$ maps from $B^n$ into ${\R}^m$ which take values in $N^k$ almost everywhere.
The critical points $u$ in $W^{1,2}(B^n,N^k)$ of the Dirichlet energy $E(u)=\int_{B^n}|\nabla u|^2$
for all perturbations of the form $\pi_{N}(u+t\phi)$, where $\phi$ is an arbitrary smooth compactly
supported map from $B^n$ into ${\R}^m$, are the harmonic maps from $B^n$ into $N^k$. They are 
the maps in $W^{1,2}(B^n,N^k)$ which solve the following Euler-Lagrange equation
(in distributional sense)
\be
\label{I.5}
-\Delta u= A(u)(\nabla u,\nabla u)=\sum_{l=1}^nA(u)(\p_{x_l}u,\p_{x_l}u)=0\quad ,
\ee
where $A(u)$ is the second fundamental form at $u(x)$ for the submanifold $N^k$ in ${\R}^m$.
For instance for $k=m-1$ when $N^{m-1}$ is an oriented codimension 1 submanifold, if we denote
by $n(y)$ the Gauss map of $N^{m-1}$ at $y$, the unit perpendicular vectorfields which generates
the orientation of $N^{m-1}$, (\ref{I.5}) becomes
\be
\label{I.6}
-\Delta u= n\, \nabla n\cdot\nabla u\quad ,
\ee
where we keep denoting $n$ the composition $n\circ u$. In order to try to extend the above described results
established for $W^{1,2}$ solutions to (\ref{I.1}) to solutions of (\ref{I.6}), or even more generally
to solutions to (\ref{I.5}), it is then natural to look for conservation laws (divergence free
quantities) generalizing (\ref{I.2}). Is for instance the non-linearity $n\, \nabla n\cdot\nabla u$ in the
right-hand-side of (\ref{I.6}) (or even (\ref{I.5}) in the local Hardy space ${\mathcal H}^1_{loc}$ ?
Do we have estimates of the form (\ref{I.4}) for general $W^{1,2}$ solutions to
(\ref{I.6}) or even (\ref{I.5})? {\bf  Can we write the equation (\ref{I.6}) or (\ref{I.5}) in divergence form ?}
Until now the answer to these questions were open and only the introduction of the indirect but beautiful technic of moving frame by F.H\'elein permitted to avoid the {\it direct conservation law approach}
for proving questions like the regularity in 2 dimension of the harmonic maps into general target 
(i.e. solutions to (\ref{I.6})) - see again \cite{Hel}.
This set of questions have motivated the following, which is one of the main result of the present paper :

\newpage

\begin{Th}
\label{th-I.1}
Let $m\in {\N}$. For every $\Omega=(\Omega_j^i)_{1\le i,j\le m}$
in $L^2(D^2,so(m)\otimes{\R}^2)$ (i.e. $\forall i,j\in\{1,\cdots ,m\}$, $\Omega_j^i\in L^{2}(D^2,{\R}^2)$
and $\Omega_j^i=-\Omega_i^j$),  every $u\in W^{1,2}(D^2,{\R}^m)$ solving
\be
\label{I.7}
-\Delta u=\Omega\cdot\nabla u
\ee
is continuous where the contracted notation in (\ref{I.7}) using coordinates stands for $\forall i=1\cdots m$ $-\Delta u^i=\sum_{j=1}^m\Omega_j^i\cdot\nabla u^j$.
\end{Th}
This theorem applies to equations (\ref{I.6}) for instance because of the following obsevation :
every derivative of $u$ solving (\ref{I.6}) is tangent to $N^{m-1}$ and is therefore perpendicular
to $n$. Thus $\sum_{j=1}^{m}n_j \nabla u^j=0$ and we can rewrite (\ref{I.6}) in the form :
\be
\label{I.8}
-\Delta u^i=\sum_{j=1}^{m}\lf[n^i\nabla n_j-n_j\nabla n^i\rg]\cdot\nabla u^j
\ee
Taking now $\Omega_j^i:=n^i\nabla n_j-n_j\nabla n^i$ we can apply theorem~\ref{th-I.1} to get the continuity of $u$.
This way of rewriting the equation has to be compared with the particular case (\ref{I.3}) except that in the general
case there is no reason for $\Omega_j^i:=n^i\nabla n_j-n_j\nabla n^i$ to be divergence free.
One of the main observation of the present work is that what is important in (\ref{I.3}) is not
the divergence free structure of $n^i\nabla n_j-n_j\nabla n^i$, valid in the particular case of the round
sphere and which disapear as soon soon as one perturbs the metric of the target, but it is the anti-symmetry of this quantity which is much more robust and which is the key point for the regularity of 
solution to (\ref{I.6}). This is a new compensation phenomenon that we discovered which goes beyond the curl-grad structures although it is strongly linked to it as we will explain in the paper.
In fact we observed that not only solutions to (\ref{I.6}), not only solutions to (\ref{I.5}) but every critical point of any elliptic conformally invariant lagrangian in 2 dimension
can be written in the form ({\ref{I.7}) and the regularity result obtained in theorem~\ref{th-I.1} can be
applied to them.
Precisely we have
\begin{Th}
\label{th-I.2}
Let $N^k$ be a $C^2$ submanifold of ${\R}^m$ ($k$ and $m$ being arbitrary integer
satisfying $1\le k\le m$). Let $\om$ be a $C^1$ $2-$form on ${N}^k$ such that the $L^\infty$ norm of $d\om$ is bounded on $N^k$. Then every critical point in $W^{1,2}(D^2,N^k)$ of the Lagrangian
\be
\label{I.9}
F(u)=\int_{D^2}\lf[|\nabla u|^2+\om(u)(\p_xu,\p_yu)\rg]\ dx\wedge dy
\ee
satisfies an equation of the form (\ref{I.7}) for some $\Om$ in $L^2(D^2,so(m)\otimes{\R}^2)$ and is therefore continuous.
\end{Th}
Critical points of $F$ which are conformal are immersed discs in $N^k$ whose mean curvature in $N^k$ at $u$ is given by $|\nabla u|^{-2}d\om(u)(\cdot,\p_xu,\p_yu)$. This is the so called prescribed mean curvature equation in a manifold $N^k$. It is not difficult to see that Lagrangian of the form (\ref{I.9})
are conformally invariant. Conversely,  It was  proved in \cite{Gr1} that every conformally invariant
elliptic Lagrangian, satisfying some "natural conditions", generates an Euler-Lagrange equation 
corresponding to a prescribed mean curvature equation in a manifold.
A particular case of interest is the case $k=m=3$ and ${N}^3={\R}^3$. Denote 
$2d\om=H(z) dz_1dz_2dz_3$ the Euler-Lagrange equation to $F$ in that case is
\be
\label{I.10}
-\Delta u=-2H(u)\,\p_xu\wedge\p_yu\quad .
\ee
There has been several attempts to prove the continuity of solutions to (\ref{I.10}) under several assumptions on $H$ like $\|H\|_{W^{1,\infty}({\R}^3)}<+\infty$
(see for instance  \cite{Hei1}, \cite{Hei2}, \cite{Gr2}, \cite{Bet1}, \cite{BeG1}, \cite{BeG2}). It was conjectured by E.Heinz, see \cite{Hei3}, that the weakest possible assumption $\|H\|_{L^\infty({\R}^3)}<+\infty$ should suffices
to ensure the continuity of $W^{1,2}$ solutions to (\ref{I.10}).
Denote $\nabla^\perp:=(-\p_y,\p_x)$ and introducing
\be
\label{I.11}
\Om:=H(u)\lf(
\begin{array}{ccc}
0& \nabla^\perp u^3&-\nabla^\perp u^2\\[5mm]
-\nabla^\perp u^3&0&\nabla^\perp u^1\\[5mm]
\nabla^\perp u^2&-\nabla^\perp u^1&0
\end{array}
\rg)
\ee
 equation (\ref{I.10}) becomes of the form
 \[
 -\Delta u=\Om\cdot\nabla u\quad,
 \]
 where $\Om\in L^2(D^2,so(3)\otimes{\R}^2)$. We can then apply theorem~\ref{th-I.1}
 to  (\ref{I.10}) and we have then proved Heinz's conjecture on prescribed mean curvature equations.

Theorem~\ref{th-I.1} is based on the discovery of conservation laws generalizing (\ref{I.2}).
 Denoting  ${M}_m({\R})$ the space of square $m\times m$ real matrices, we have :
 \begin{Th}
\label{th-I.3}
Let $m\in {\N}$. Let $\Omega=(\Omega_j^i)_{1\le i,j\le m}$
in $L^2(B^n,so(m)\otimes\wedge^1{\R}^n)$ and let $A\in L^\infty(B^n,{ M}_m({\R}))\cap
W^{1,2}$ and $B\in W^{1,2}(B^n,{M}_m({\R})\otimes\wedge^2{\R}^n)$ satisfying
\be
\label{I.12z}
d_\Om A:= d A- A \Om=-d^\ast B
\ee
(where explicitely (\ref{I.12}) means $\forall i,j\in\{1\cdots m\}$ $d A^i_j-\sum_{k=1}^mA^i_k\Om^k_j=-d^\ast B^i_j$)
Then every solution to (\ref{I.7}) on $B^n$ satisfies the following conservation law
\be
\label{I.13z}
d\lf(\ast A\ du+(-1)^{n-1}(\ast B)\wedge du\rg)=0\quad .
\ee
\end{Th}
For $n=2$, using different notations, the theorem says that given $\Omega=(\Omega_j^i)_{1\le i,j\le m}$
in $L^2(D^2,so(m)\otimes{\R}^2)$, $A\in L^\infty(D^2,{ M}_m({\R}))\cap
W^{1,2}$ and $B\in W^{1,2}(D^2,{M}_m({\R}))$ satisfying
\be
\label{I.12}
\nabla_\Om A:=\nabla A- A \Om=\nabla^\perp B
\ee
Then every solution to (\ref{I.7}) satisfies the following conservation law
\be
\label{I.13}
div(A\nabla u+B\nabla^\perp u)=0\quad .
\ee
For instance going back to the symmetric situation of harmonic maps into $S^{m-1}$ we take 
$A$ and $B$ satisfying
\[
\lf\{
\begin{array}{l}
\ds A=id_m=(\delta^j_i)_{1\le i,j\le m}\\[5mm]
\ds\nabla^\perp B^i_j=u^i\nabla u_j-u_j\nabla u^i\quad. 
\end{array}
\rg.
\]
Using now the fact that $div( B^i_j\nabla^\perp u^j)=\nabla B^i_j\cdot\nabla^\perp u^j=
-\nabla^\perp B^i_j\cdot\nabla u^j$ the harmonic map equation into $S^{m-1}$, with these notations,
is equivalent to (\ref{I.13}). We then have included the classical conservation law (\ref{I.2})
into the larger set of conservation laws of the form (\ref{I.13}). The question remains of finding
$A$ and $B$ satisfying (\ref{I.12}). We shall prove the following local existence result 
\begin{Th}
\label{th-I.4}
There exists $\ep(m)>0$ and $C(m)$ such that, for every  $\Omega=(\Omega_j^i)_{1\le i,j\le m}$
in $L^2(D^2,so(m)\otimes{\R}^2)$ satisfying
\be
\label{I.14}
\int_{D^2}|\Om|^2\le\ep_m\quad,
\ee
there exists $A\in L^\infty(D^2,Gl_m({\R}))\cap
W^{1,2}$ and $B\in W^{1,2}(D^2,{ M}_m({\R}))$ satisfying 
\begin{itemize}
\item[i)] 
\be
\label{I.15}
\int_{D^2}|\nabla A|^2+|\nabla A^{-1}|^2 +\|dist(A,SO(n))\|_\infty^2\le C(n)\int_{D^2}|\Om|^2\quad,
\ee
\item[ii)]
\be
\label{I.16}
\int_{D^2}|\nabla B|^2\le C(n)\int_{D^2}|\Om|^2\quad,
\ee
\item[iii)]
\be
\label{I.17}
\nabla_\Om A:=\nabla A- A \Om=\nabla^\perp B\quad .
\ee
\end{itemize}
\end{Th}
A corresponding local existence result in higher dimension is still an open problem.
If the harmonic map we are considering is stationary (see \cite{Hel}), $\Om$ is in the Morrey space
 given by
 \be
 \label{I.27}
\|\Om\|_{M^1_2}=sup_{x,r}\frac{1}{r^{n-2}}\int_{B_r(x)}|\nabla \Om|<+\infty\quad.
 \ee
Then, under the assumption that $\|\Om\|_{M^1_2}$ is below some positive constant depending only on $n$ and $m$, the elliptic linear system (\ref{I.12z}) becomes critical and the  existence of $A$ and $B$ solving (\ref{I.12z}) should be looked for in the space
${M^2_2}$ ( following the search of a Coulomb gauge in Morrey spaces introduced in \cite{MeR}, one has a replacement of lemma~\ref{lm-A3} in higher dimension).

Local existence of conservation law (\ref{I.13}) for stationary harmonic maps permits to extend to general $C^2$ targets 
 the partial regularity of  L.C. Evans \cite{Ev} for harmonic maps into spheres  following the same strategy that Evans
 introduced. 

 Using conservation laws (\ref{I.13}) we can  prove the following result.
\begin{Th}
\label{th-I.4b}
Let $\Om_n\in L^2(D^2,so(m)\otimes\wedge^1{\R}^2)$ such that $\Om_n$ weakly converges
in $L^2$ to some $\Om$. Let $f_n$ be a sequence in $H^{-1}(D^2,{\R}^m)$ which converges
to 0 in $H^{-1}$ and $u_n$ be a bounded sequence in $W^{1,2}(D^2,{\R}^m)$ solving
\be
\label{I.171}
-\Delta u_n=\Omega_n\cdot\nabla u_n+f_n\quad\quad\mbox{ in } D^2
\ee  
Then, there exists a subsequence $u_{n'}$ of $u_n$ which weakly converges in $W^{1,2}$ to a solution of 
(\ref{I.7}).
\end{Th}
Passage to the limit in the equation in 2 dimension for the prescribed mean curvature equation
or for the harmonic map equation was established in \cite{Bet2} using involved technics. A much simpler proof
using moving frames was then given in \cite{FMS}. In both proofs a Lipschitz bound on the prescribed mean curvature 
was required. This is no more the case in theorem~\ref{th-I.4b} where only an $L^\infty$ bound on the 
prescribed mean curvature is needed.

Following similar ideas theorem~\ref{th-I.1} can be extended in it's spirit to first order ellitic complex valued PDE
\begin{Th}
\label{th-I.5}
Let $m\in {\N}$ and $k\in{\N}$.  Let $\Omega=(\Omega_j^i)_{1\le i,j\le m}$
in $L^2(D^2,so(m)\otimes{\C}\otimes\wedge^1{\R}^2)$ and let $\al\in L^2(D^2,M_{m,k}({\C}))$ solving
\be
\label{I.18}
\frac{\p\al}{\p\ov{z}}=\Om\,\al\quad,
\ee
then there exists $P\in W^{1,2}(D^2, SO_m({\C}))$ and $\beta\in C^\infty(D^2,M_{m,k}({\C})$
such that $\al=P\beta$, where $SO_m({\C})$ is the group of invertible matrices in $Gl_m({\C})$
satisfying $P^tP=id_{m}$.
\end{Th}

\noindent{ \bf Conservation laws and moving frames.}

Existence of global conservation laws can be obtained in the same spirit by the mean of moving frames.
Considering a map $u$  in $W^{1,2}(D^2,N^2)$ where  $N^2$ is a closed $C^2$ oriented 2-dimensional submanifold of ${\R}^m$, then there exists a map $e_1$ in $W^{1,2}(D^2,S^{m-1})$ such that $e_1(x)\in T_{u(x)}N^2$
for almost every $x$ in $D^2$. Moreover, denoting $e_2(x)$ the unit vector perpendicular to $e_1$ such that $e_1\wedge e_2$ is the unit 2-vector giving the oriented tangent plane $T_{u(x)}N^2$, we can choose $e_1$ such that $div((e_2,\nabla e_1))=0$ on $D^2$ (see \cite{Hel}) where $(\cdot,\cdot)$ 
denotes the scalar product in ${\R}^m$ that we also sometime simply denote $\cdot$. Such a pair $(e_1,e_2)$
is called a Coulomb moving frame associated to $u$. We have then the following conservation law
\begin{Th}
\label{th-I.6}
Let $u$ be a $W^{1,2}$ harmonic map from $D^2$ into $N^2$, let $(e_1,e_2)$ be a Coulomb moving frame associated to $u$, let $a$ be the function
solving
\be
\label{I.19}
\lf\{
\begin{array}{l}
\ds -\Delta a=(\nabla^\perp e_1,\nabla e_2)=\frac{\p e_1}{\p x}\cdot\frac{\p e_2}{\p y}-
\frac{\p e_1}{\p y}\cdot\frac{\p e_2}{\p x}\quad\mbox{ in }D^2\\[5mm]
\ds a=0\quad\quad\mbox{ on }\p D^2
\end{array}
\rg.
\ee
then the following conservation law holds
\be
\label{I.20}
\lf\{
\begin{array}{l}
\ds div\lf(ch a\, (\nabla u,e_1)+sh a\, (\nabla^\perp u,e_2)\rg)=0\quad ,\\[5mm]
\ds div\lf(ch a\, (\nabla u,e_2)-sh a\, (\nabla^\perp u,e_1)\rg)=0\quad .
\end{array}
\rg.
\ee
Moreover the following estimate holds
\be
\label{I.21}
\int_{D^2_{1/2}}|\nabla^2 u|\le C \exp\lf[\frac{1}{4\pi}\int_{D^2}|\nabla e|^2\rg]\, \lf(\|\nabla e\|_{L^2(D^2)}+1\rg)\ \|\nabla u\|_{L^2(D^2)}\quad,
\ee
where $|\nabla e|^2:=|\nabla e_1|^2+|\nabla e_2|^2$.
\end{Th}
Observe that, because of Wente's Lemma~\ref{lm-A1} that we recall in the appendix, the solution $a$
of (\ref{I.19}) is bounded in $L^\infty$.
When $N^2$ is not diffeomorphic to $S^2$ one can estimate $\int_{D^2}|\nabla e|^2$ in terms
of $\|u\|_{W^{1,2}}$. This is no more the case for $N^2=S^2$ : one can find sequences of $u_n$,
harmonic from $D^2$ into $S^2$ with uniformly bounded $W^{1,2}$ norm but for which, however,
every sequence of Coulomb moving frame is not bounded in $W^{1,2}$. Nevertheless we
still believe that the following holds true :

{\bf Conjecture : }{\it For every $k\le m$, for every $n\in {\N}$ for every $N^k$, k-dimensional closed submanifold of ${\R}^m$, and for every $C>0$ there exists $\delta(C,n,N^k)>0$ such that if $u$ is a $W^{1,2}$ harmonic map from $B^n_2(0)$
into $N^k$ satisfying
\be
\label{I.22}
\int_{B^n_2(0)}|\nabla u|^2\le \delta(C,n,N^k)
\ee
then
\be
\label{I.23}
\int_{B^n_1(0)}|\nabla^2 u|\le C
\ee}
Such an estimate would, in particular, permit to extend the quantization result of \cite{LiR} to general targets.

Finally, in general dimension, the following conservation laws generalizing (\ref{I.20})
should play an important role in the theory of harmonic maps : 
\begin{Th}
\label{th-I.7}
Let $u$ be a $W^{1,2}$ harmonic map from $B^n$ into $N^k $, a closed oriented $C^2$ $k-$dimensional
submanifold of ${\R}^m$. Let $(e_1,\cdots, e_k)$ be a Coulomb moving frame associated
to $u$  ( the map $x\rightarrow (e_1,\cdots,e_k)$ is in $W^{1,2}$, for almost every $x$
$(e_1(x),\cdots,e_k(x))$ is an orthonormal basis of $T_{u(x)}N^k$ and $\forall i,j\in\{1\cdots k\}$
$d^\ast(e_i,de_j)=0$.) Denote $\Om=(\Om_i^j)\in L^2(B^n,so(k)\otimes\wedge^1{\R}^n)$ the connection
given by
\[
\Om_{i}^j:=(e_j,de_i)\quad .
\]
Let $\Phi\in L^4(B^n,M_k({\R}))\cap W^{1,2}$ and $\Psi\in L^4(B^n,M_k({\R})\otimes\wedge^2{\R}^n)
\cap W^{1,2}$ solving
the linear equation
\be
\label{I.24}
d_\Om\Phi+d^\ast_\Om\Psi=0\quad,
\ee
where $(d_\Om\Phi)^j_i:=d\Phi^j_i+\Phi_i^k\wedge\Om_k^j$ and $d^\ast_\Om$ is the adjoint
of $d_\Om$ given by $(d^\ast_\Om\Psi)^j_i:=d^\ast\Psi^j_i+\ast(\ast\Psi_i^k\wedge\Om_k^j)$.
Then the following conservation laws are satisfied
\be
\label{I.25}
d\lf(\ast\Phi (du,e)+(-1)^{n-1}(\ast\Psi)\wedge (du,e)\rg)=0\quad .
\ee
where $(du,e)$ is the element in $L^2(B^n,{\R}^k\otimes\wedge^1{\R}^n)$ given by
$\{(du,e_j)\}_{j=1\cdots k}$ and $\Phi (du,e)$ and $\ast\Psi\wedge (du,e)$ denote respectively
the elements  in $L^\frac{4}{3}(B^n,{\R}^k\otimes\wedge^1{\R}^n)$ and in 
$L^\frac{4}{3}(B^n,{\R}^k\otimes\wedge^{n-1}{\R}^n)$ given by $\Phi_i^j(du,e_j)$ and 
$\ast\Psi_i^j\wedge(du,e_j)$.
\end{Th}
 Observe that (\ref{I.25}) generalizes (\ref{I.20}) to general dimension by taking
 \be
 \label{I.26}
\ds \Phi=ch a\  \lf(\begin{array}{cc}
 1& 0\\[5mm]
 0& 1
 \end{array}
 \rg)\quad\mbox{ and }\quad 
\Psi=sh a\  \lf(\begin{array}{cc}
 0& -1\\[5mm]
 1& 0
 \end{array}
 \rg)\ dx\wedge dy
 \ee
 Existence of Coulomb moving frames is discussed in \cite{Hel} and is proved under the assumption that $N^k$ is sufficiently regular and modulo some isometric embeddings in a submanifold diffeomorphic 
 to a torus. Again here, under the assumption that $\|\Om\|_{M^1_2}$ is below some positive constant depending only on $n$ and $k$, the elliptic linear system (\ref{I.24}) becomes critical and the  existence of $\Phi$ and $\Psi$ solving (\ref{I.24}) should be looked for in the space
${M^2_2}$ which embbeds in $L^4$ in every dimension. 

The paper is organised as follows : in section 2 we prove theorem~\ref{th-I.1} to theorem~\ref{th-I.5}.
In section 3 we prove theorem~\ref{th-I.4b}. In section 4 we prove theorem~\ref{th-I.6} and theorem~\ref{th-I.7}.
In the appendix we recall Wente's result and establish several lemmas used in sections 2 and 3.

\noindent{\bf Acknowledgments :} {\it This work was carried out while the author was visiting the
Universit\'e de Bretagne Occidentale at Brest. The author would have like to thank the mathematic
department of the UBO for it's hospitality.}

\section{Proof of theorems~\ref{th-I.1}...\ref{th-I.5}.}
\reset
\subsection{ Proof of theorem~\ref{th-I.4}.}
Let $\ep(m)>0$ given by lemma~\ref{lm-A3}. Let $\Om\in L^2(D^2,so(m)\otimes\wedge^1{\R}^2)$
satisfying
\be
\label{II.1}
\int_{D^2}|\Om|^2\le\ep(m)\quad .
\ee
Let then $P\in W^{1,2}(D^2,SO(m))$ and $\xi\in W^{1,2}(D^2,so(m))$ given by lemma~\ref{lm-A3}
satisfying
\be
\label{II.2}
\nabla^\perp\xi=P^{-1}\nabla P+P^{-1}\Om P\quad\mbox{ in }D^2\quad,
\ee
with $\xi=0$ on $\p D^2$ and such that
\be
\label{II.3}
\|\xi\|_{W^{1,2}}+\|P\|_{W^{1,2}}+\|P^{-1}\|_{W^{1,2}}\le C(m)\ \|\Om\|_{L^2}\quad.
\ee
We look for $A$ and $B$ solving (\ref{I.12}) and introducing $\ti{A}:=AP$, it means that we are looking for
$\ti{A}$ and $B$ solving
\be
\label{II.4}
\nabla\ti{A}-\nabla^\perp B\  P=\ti{A}\nabla^\perp\xi\quad .
\ee
First we aim to solve the following system
\be
\label{II.5}
\lf\{
\begin{array}{l}
\ds \Delta \hat{A}=\nabla\hat{A}\cdot\nabla^\perp\xi-\nabla^\perp B\cdot\nabla P\quad\mbox{ in }D^2\quad,\\[5mm]
\ds \Delta B=\nabla^\perp\hat{A}\cdot\nabla P^{-1}-div(\hat{A}\nabla\xi P^{-1})-div(\nabla\xi P^{-1})\quad\mbox{ in }D^2\quad,\\[5mm]
\ds \frac{\p\hat{A}}{\p\nu}=0\quad\mbox{ and }\quad B=0\quad\mbox{ on }\p D^2\\[5mm]
\ds\mbox{ and }\int_{D^2}\hat{A}=id_{m}
\end{array}
\rg.
\ee
for $\hat{A}\in M_m({\R})$ and $B\in{M_m({\R})}$. Observing that the right-hand-side of the first
equation of (\ref{II.5}) is made of jacobians : $(\nabla\hat{A}\cdot\nabla^\perp\xi)^j_i=
\p_y\hat{A}_i^k\, \p_x\xi_k^j-\p_x\hat{A}_i^k\,\p_y\xi_k^j$ and $-(\nabla^\perp B\cdot\nabla P)_i^j=
\p_xB_i^k\,\p_y P^j_k-\p_yB_i^k\,\p_xP^j_k$, using lemma~\ref{lm-A1}, standard elliptic estimates
and the fact that $P\in SO(m)$, we have the a-priori estimates
\be
\label{II.6}
\ds\|\hat{A}\|_{W^{1,2}}+\|\hat{A}\|_{L^\infty}\le C\ \|\xi\|_{W^{1,2}}\ \|\hat{A}\|_{W^{1,2}}
+ C\ \|P\|_{W^{1,2}}\ \|B\|_{W^{1,2}}
\ee
\be
\label{II.6q}
\ds\|B\|_{W^{1,2}}\le C\ \|P^{-1}\|_{W^{1,2}}\ \|\hat{A}\|_{W^{1,2}}
+ C\ \|\xi\|_{W^{1,2}}\ \|\hat{A}\|_{L^\infty}+ C \|\xi\|_{W^{1,2}}
\ee
Thus for $\|\Om\|_{L^2}$ small enough, using a standard fixed point argument, we obtain the existence of
$\hat{A}$ and $B$ satisfying (\ref{II.5}) and
\be
\label{II.7}
\|\hat{A}\|_{W^{1,2}}+\|\hat{A}\|_{L^\infty}+\|B\|_{W^{1,2}}\le C\ \|\Om\|_{L^2}
\ee
(Observe that, using the result of \cite{CLMS}, we even have $\|\hat{A}\|_{W^{2,1}}\le C\ \|\Om\|_{L^2}$).
Let now $\ti{A}:=\hat{A}+id_m$. From the first equation of (\ref{II.5}) we obtain the existence of $C$ in $W^{1,2}$ satisfying
\be
\label{II.8}
\nabla\ti{A}-\ti{A}\,\nabla^\perp\xi-\nabla^\perp B=\nabla^\perp C\quad.
\ee
Moreover, $C$ satisfies
\be
\label{II.9}
\lf\{
\begin{array}{l}
\ds div(\nabla C\ P^{-1})=0\quad\mbox{ in }D^2\\[5mm]
\ds C=0\quad\mbox{ on }\p D^2\quad .
\end{array}
\rg.
\ee
Using now lemma~\ref{lm-A2}, we obtain that $C$ is identically zero and $A:=\ti{A} P^{-1}$ and $B$
satisfy (\ref{I.15}), (\ref{I.16}) and (\ref{I.17}). Theorem~\ref{th-I.4} is then proved.

\subsection{Proof of theorem~\ref{th-I.3}.}

Theorem~\ref{th-I.3} follows from a direct computation.

\subsection{Proof of theorem~\ref{th-I.1}.}

Since the desired result is local (continuity of $u$), we can always assume that $\int_{D^2}|\Om|^2\le\ep(m)$
where $\ep(m)$ is given by theorem~\ref{th-I.4}. Moreover, let $A$ and $B$ given by theorem~\ref{th-I.4}.
From theorem~\ref{th-I.4} they solve the following system
\be
\label{II.11}
\lf\{
\begin{array}{l}
\ds div(A\,\nabla u)=-\nabla B\cdot\nabla^\perp u\quad,\\[5mm]
\ds curl(A\,\nabla u)=\nabla^\perp A\cdot\nabla u
\end{array}
\rg.
\ee
Using standard  elliptic estimates, we get the existence of $E$ and
$D$ in $W^{1,2}(D^2)$ such that
\be
\label{II.12}
A\,\nabla u=\nabla ^\perp E+\nabla D\quad .
\ee
Moreover, using the jacobian structure of the right-hand-sides of the equations in (\ref{II.11}), the results in \cite{CLMS} imply that $E$ and $D$ are in $W^{2,1}$ on the disk of half radius $D^2_{1/2}$.
Therefore $\nabla u=A^{-1}\nabla^\perp E+A^{-1}\nabla D$ is in $W^{1,1}$ on this disk.
Using the embedding of $W^{2,1}$ into $C^0$ in 2 dimension we get the desired result and theorem~\ref{th-I.1} is proved.

\subsection{Proof of theorem~\ref{th-I.2}.}
Let $N^k$ be a $k-$dimensional submanifold of ${\R}^m$. Let $\pi_N$ be the orthogonal projection on $N^k$ defined in a small tubular neighborhood of $N^k$. Let $\om$ be a 2-form on $N^k$ 
and let $\ti{\om}$ be the pull back of $\om$ by $\pi_N$ in this small tubular neighborhood : $\ti{\om}:=\pi_N^\ast\om$.
Following \cite{Hel} chapter 4, critical points of (\ref{I.9}) in $W^{1,2}(D^2,N^k)$ satisfy the following
Euler-Lagrange equation 
\be
\label{II.13}
\Delta u^i+A^i(u)_{j,l}\nabla u^j\cdot\nabla u^l+\la(u)^i_{j,l}\p_x u^j\p_y u^l=0
\ee
where $\la(u)^i_{j,l}:=d\ti{\om}_{u}(\ep_i,\ep_j,\ep_l)$ where $(\ep_l)_{l=1\cdots m}$
is the canonical basis of ${\R}^m$. Thus, in particular, $\la(u)^i_{j,l}=-\la(u)^j_{i,l}$.
Since $(A^j_{i,l})_{j=1\cdots m}=A(\ep_i,\ep_l)$ is
perpendicular to $T_uN^k$ for every $i$ and $l$, we have that
\be
\label{II.14}
\forall i,l\in \{1\cdots m\}\quad\quad \sum_j A^j_{i,l}\nabla u^j=0\quad .
\ee
Thus finally (\ref{II.13}) becomes
\be
\label{II.15}
\begin{array}{l}
\ds\Delta u^i +\lf[A^i(u)_{j,l}-A^j(u)_{i,l}\rg] \nabla u^l\cdot\nabla u^j\\[5mm]
\ds\quad\quad+\frac{1}{4}\lf[\la(u)^i_{j,l}-\la(u)^j_{i,l}\rg]\nabla^\perp u^l\cdot\nabla u^j=0\quad.
\end{array}
\ee
Introducing $\Om:=(\Om^j_i)_{i,j}$ where
\be
\label{II.16}
\Om^i_j:=\lf[A^i(u)_{j,l}-A^j(u)_{i,l}\rg]\nabla u^l+\frac{1}{4}\lf[\la(u)^i_{j,l}-\la(u)^j_{i,l}\rg]\nabla^\perp u^l
\ee
we have succeeded in writing $W^{1,2}$ critical points of lagrangian of the form (\ref{I.9}) as solutions
to PDE of the form (\ref{I.7}) for some $\Om\in L^2(D^2,so(m)\otimes\wedge^1{\R}^2)$. 
Theorem~\ref{th-I.2} is then proved.

\section{Conservation laws and passage to the limit in PDEs.}
\reset

The goal of this section is to prove theorem~\ref{th-I.4b}.

Let $\Om_n\in L^2(D^2,so(m)\otimes\wedge^1{\R}^2)$ converging weakly to some $\Om$ and $f_n$ and $u_n$
respectively converging to zero in $H^{-1}(D^2,{\R}^m)$ and bounded in $W^{1,2}(D^2,{\R}^m)$. We can allways assume
that $u_n$ converges weakly to some $u$ in $W^{1,2}(D^2,{\R}^m)$. Let $\la<1$ and let $\ep(m)$ given by theorem~\ref{th-I.4}. 
To every $x$ in $B_\la^2(0)$ we assign $r_{x,n}\le 1-|x|$ such that $\int_{B_{r_{x,n}}(x)}|\Om_n|^2=\ep(m)$ or $r_{x,n}=1-|x|$
in case $\int_{B_{1-|x|}(x)}|\Om_n|^2<\ep(m)$. $\{B_{r_x}(x)\}$ for every $x$ in $B_\la^2(0)$ realizes of course a covering
of $B_{\la}^2(0)$. We extract a Vitali covering from it which ensures that every point in  $B_{\la}^2(0)$ is covered by
a number of ball bounded by a universal number. Since $\int_{D^2}|\Om_n|^2$ is uniformly bounded, the number of balls in each
such a Vitali covering for each $n$ is also uniformly bounded and, modulo extraction of a subsequence, we can assume
that it is fixed and equal to $N$ independent of $n$. Let $\{B_{r_{i,n}}(x_{i,n})\}_{i=1\cdots N}$ be this covering.
Modulo extraction of a subsequence we can allways assume that each sequence $x_{i,n}$ converges in $\ov{B}_\la^2(0)$
to a limit $x_i$ and that each sequence $r_{i,n}$ converges to a non negative number $r_{i}$ (which could be zero of course).
We claim that equation (\ref{I.7}) is satisfied on each $B_{r_i}(x_i)$. Let $A_{i,n}$ and $B_{i,n}$ given by theorem~\ref{th-I.4}
in $B_{r_{i,n}}(x_{i,n})$ for $\Om_n$. We then have   
\be
\label{II.16b}
div(A_{i,n}\nabla u_n+B_{i,n}\nabla^\perp u_n)=-A_{i,n}\ f_n\quad\quad\mbox{ in }B_{r_{i,n}}(x_{i,n})\quad.
\ee
where $A_{i,n}$ and $B_{i,n}$ satisfy
\be
\label{II.16c}
\nabla A_{i,n}- A_{i,n} \Om_{i,n}=\nabla^\perp B_{i,n}
\ee
We can extract a subsequence such that each of the couples $(A_{i,n},B_{i,n})$ weakly converge in $W^{1,2}$
to some limit $(A_i,B_i)$
in every $B_{ r_i}(x_i)$. Because of the weak convergences in $W^{1,2}$ we have 
strong convergences in $L^2$ and then we have that
\be
\label{II.17b}
A_{i,n}\nabla u_n+B_{i,n}\nabla^\perp u_n\longrightarrow A_i \nabla u+B_i\nabla^\perp u\quad\quad\mbox{ in }{\mathcal D}'
\ee
\be
\label{II.18b}
\nabla A_{i,n}- A_{i,n} \Om_{i,n}-\nabla^\perp B_{i,n}\longrightarrow \nabla A_{i}- A_{i} \Om-\nabla^\perp B_{i}\quad\quad\mbox{ in }{\mathcal D}'
\ee
and
\be
\label{II.19b}
-A_{i,n}\ f_n\rightarrow 0\quad\quad\mbox{ in }{\mathcal D}'
\ee
Combining then (\ref{II.16b})...(\ref{II.19b}) we obtain that
\be
\label{II.20b}
div(A_{i}\nabla u+B_{i}\nabla^\perp u)=0\quad\mbox{ in }B_{r_i}(x_i)\quad,
\ee
and that
\be
\label{II.21b}
\nabla A_i-A_i\Om=\nabla^\perp B_i\quad\mbox{ in }B_{r_i}(x_i)\quad.
\ee
Combining (\ref{II.20b}) and (\ref{II.21b}) we then have that
\be
\label{II.22b}
A_i\lf[\Delta u+\Om\cdot \nabla u\rg]=0\quad\mbox{ in }B_{r_i}(x_i)\quad.
\ee
From (\ref{I.15}), because of the pointwise convergence of $A_{i,n}$, we get the invertibility
of $A_i$ and (\ref{II.22b}) implies that
\be
\label{II.23}
\Delta u+\Om\cdot \nabla u=0\quad\mbox{ in }B_{r_i}(x_i)
\ee
It is clear that every point in $B_{\la}^2(0)$ is in the closure of the union of the $B_{r_i}(x_i)$. Let $x$ be a point
which is none of the $B_{r_i}(x_i)$. It seats then on the circle, boundary of one of the $B_{r_i}(x_i)$.
For convexity reason, it has to seat at the boundary of at least 2 different circles. 2 different circles can intersect
at only finitely many points (0,1 or 2 points), since there are finitely many circles, only finitely many points in $B_{\la}^2(0)$
can be outside the union of the $B_{r_i}(x_i)$.  Thus the distribution $\Delta u+\Om\cdot \nabla u$ is supported at at mostly
 finitely many points. Since $\Delta u+\Om\cdot \nabla u\in H^{-1}+L^{1}$ it is identically zero on $B_{\la}^2(0)$.
Since this holds for every $\la<1$ we have proved theorem~\ref{th-I.4b}.

\section{Conservation laws and moving frames.}
\reset
\subsection{Proof of theorem~\ref{th-I.6}.}

First (\ref{I.20}) is the result of a direct computation, granting the fact that $(\Delta u,e)=0$.
It remains then to prove (\ref{I.21}). We rewrite (\ref{I.20}) in the form (using ${\Z}_2$ indexation
\be
\label{III.1}
\lf\{
\begin{array}{l}
\ds div (ch\,a (\nabla u,e_i))=(-1)^i(\nabla(sh\, a\  e_{i+1}),\nabla^\perp u)\quad,\\[5mm]
\ds curl(ch\,a (\nabla u, e_i))=(\nabla u,\nabla^\perp(ch\,a\  e_i))\quad .
\end{array}
\rg.
\ee
Using then lemma~\ref{lm-A1}, and standard elliptic estimate, there exist 
$E\in W^{2,1}(D^2_{1/2},{\R}^2)$ and $D\in W^{2,1}(D^2_{1/2},{\R}^2)$ such that
\be
\label{III.2}
ch\,a (\nabla u, e_i)=\nabla E_i+\nabla^\perp D_i\quad ,
\ee
moreover
\be
\label{III.3}
\begin{array}{l}
\ds\|E\|_{W^{2,1}(D^2_{1/2})}+\|D\|_{W^{2,1}(D^2_{1/2})}\le C\|\nabla(sh\, a\  e)\|_{L^2}\ \|\nabla u\|_{L^2}\\[5mm]
\quad\quad\ds +\|\nabla(ch\, a\  e)\|_{L^2}\ \|\nabla u\|_{L^2}+\|ch\,a (\nabla u, e)\|_{L^2}\quad,
\end{array} 
\ee
Thus we have
\be
\label{III.4}
\begin{array}{l}
\ds\|E\|_{W^{2,1}(D^2_{1/2})}+\|D\|_{W^{2,1}(D^2_{1/2})}\le C\, e^{\|a\|_\infty}\ \|\nabla e\|_{L^2}\ \|\nabla u\|_{L^2}\\[5mm]
 \quad\quad+e^{\|a\|_\infty}\ \|\nabla a\|_{L^2}\ \|\nabla u\|_{L^2}+e^{\|a\|_\infty}\ \|\nabla u\|_{L^2}
 \end{array}
 \ee
Applying lemma~\ref{lm-A1}  (with the optimal constant given in \cite{Hel}) to equation (\ref{I.19}) we have 
\be
\label{III.5}
\|a\|_{L^\infty(D^2)}\le \frac{1}{2\pi}\|\nabla e_1\|_{L^2}\ \|\nabla e_2\|_{L^2}
\ee
and
\be
\label{III.5a}
\|\nabla a\|_{L^2(D^2)}\le \frac{1}{\sqrt{2\pi}}\|\nabla e_1\|_{L^2}\ \|\nabla e_2\|_{L^2}
\ee
We rewrite (\ref{III.2}) in the form
\be
\label{III.6}
\nabla u=(ch\, a)^{-1}\lf(\nabla E_j\ e_j+\nabla^\perp D_j\ e_j\rg)\quad.
\ee
Combining then (\ref{III.4}), (\ref{III.5}) and (\ref{III.6}) we get (\ref{I.21}) and theorem~\ref{th-I.6}
is proved.
\subsection{Proof of theorem~\ref{th-I.7}.}
Theorem~\ref{th-I.7} follows from a direct computation.

\appendix
\section{Appendix}
\reset
\begin{Lma}
\label{lm-A1} \cite{We}, \cite{BrC}, \cite{Ta1}, \cite{CLMS}
Let $a$ and $b$ in $L^{1}(D^2,{\R})$ such that $\nabla a$ and $\nabla b$ are in 
$L^2(D^2)$. Let $\varphi$ be the solution of
\be
\label{A.1}
\lf\{
\begin{array}{l}
\ds \Delta \varphi=\frac{\p a}{\p x}\frac{\p b}{\p y}-\frac{\p a}{\p y}\frac{\p b}{\p x}\quad\mbox{ in }D^2\\[5mm]
\ds \varphi=0\quad \mbox{ or }\quad\frac{\p \varphi}{\p \nu}=0\quad\mbox{ on }\p D^2
\end{array}
\rg.
\ee
Then the following estimates hold
\be
\label{A.2}
\|\varphi\|_{L^\infty(D^2)}+\|\nabla\varphi\|_{L^2(D^2)}+\|\nabla^2\varphi\|_{L^1(D^2)}\le
\|\nabla a\|_{L^2(D^2)}\ \|\nabla b\|_{L^2(D^2)}\quad ,
\ee
where we choose $\int_{D^2}\varphi=0$ for the Neuman boundary data.
\end{Lma}
\begin{Lma}
\label{lm-A2}
There exists $\ep>0$ such that for every $P\in W^{1,2}(D^2,Gl(m))$ satisfying
\be
\label{A.3}
\int_{D^2}|\nabla P|^2 +|\nabla P^{-1}|^2\le \ep\quad ,
\ee
then, $C\equiv 0$ is the unique  solution  in $W^{1,2}(D^2,M_m({\R}))$ of the following problem
\be
\label{A.4}
\lf\{
\begin{array}{l}
\ds div\lf(\nabla C\ P^{-1}\rg)=0\quad\mbox{ in }D^2\\[5mm]
\ds C=0\quad\quad\mbox{ on }\p D^2
\end{array}
\rg.
\ee
\end{Lma}
{\bf Proof of lemma~\ref{lm-A2}.}
Since $C$ satisfies (\ref{A.4}), there exists $D\in W^{1,2}(D^2,M_m({\R}))$ such that 
$\nabla^\perp D=\nabla C\ P^{-1}$ and we can choose $D$ such that $\int_{D^2}D=0$.
Thus $C$ and $D$ satisfy respectively
\be
\label{A.5}
\lf\{
\begin{array}{l}
\ds\Delta C=\nabla^\perp D\cdot\nabla P\quad\mbox{ in }D^2\quad ,\\[5mm]
\ds C=0\quad\quad\mbox{ on }\p D^2\quad .
\end{array}
\rg.
\ee
and
 \be
\label{A.6}
\lf\{
\begin{array}{l}
\ds\Delta D=-\nabla^\perp C\cdot\nabla P^{-1}\quad\mbox{ in }D^2\quad ,\\[5mm]
\ds \frac{\p D}{\p \nu}=0\quad\quad\mbox{ on }\p D^2\quad .
\end{array}
\rg.
\ee
Thus, using lemma~\ref{lm-A1} and (\ref{A.3}), for $\ep$ small enough, we have
\be
\label{A.7}
\begin{array}{l}
\ds \|\nabla C\|_{L^2(D^2)}\le \frac{1}{2}\|\nabla D\|_{L^2(D^2)}\quad\mbox{ and }\\[5mm]
\ds  \|\nabla D\|_{L^2(D^2)}\le \frac{1}{2}\|\nabla C\|_{L^2(D^2)}
\end{array}
\ee
which implies that $C\equiv 0$ and $D\equiv 0$ and lemma~\ref{lm-A2}
is proved.

\begin{Lma}
\label{lm-A3}
There exist $\ep(m)>0$ and $C(m)>0$ such that for every $\Om$ in 
$L^2(D^2,so(m)\otimes\wedge^1{\R}^2)$ satisfying
\be
\label{A.8}
\int_{D^2}|\Om|^2\le\ep(m)\quad,
\ee
then there exist $\xi\in W^{1,2}(D^2,so(m))$ and $P\in W^{1,2}(D^2,SO(m))$ such that
\begin{itemize}
\item[i)]
\be
\label{A.9}
\nabla^\perp\xi=P^{-1}\nabla P+P^{-1}\Om P\quad\mbox{ in }D^2\quad,
\ee
\item[ii)]
\be
\label{A.10}
\xi=0\quad\quad\mbox{ on } \p D^2\quad ,
\ee
\item[iii)]
\be
\label{A.11}
\|\xi\|_{W^{1,2}(D^2)}+\|P\|_{W^{1,2}(D^2)}\le C(m)\ \|\Om\|_{L^2(D^2)}\quad .
\ee
\end{itemize}
\end{Lma}
In order to prove lemma~\ref{lm-A3} we follow the strategy of \cite{Uhl} and we first prove the following 
result
\begin{Lma}
\label{lm-A4}
There exist $\ep(m)>0$ and $C(m)>0$ such that for every $\Om$ in 
$W^{1,2}(D^2,so(m)\otimes\wedge^1{\R}^2)$ satisfying
\be
\label{A.12}
\int_{D^2}|\Om|^2\le\ep(m)\quad,
\ee
then there exist $\xi\in W^{2,2}(D^2,so(m))$ and $P\in W^{2,2}(D^2,SO(m))$ such that
\begin{itemize}
\item[i)]
\be
\label{A.13}
\nabla^\perp\xi=P^{-1}\nabla P+P^{-1}\Om P\quad\mbox{ in }D^2\quad,
\ee
\item[ii)]
\be
\label{A.14}
\xi=0\quad\quad\mbox{ on } \p D^2\quad,
\ee
\item[iii)]
\be
\label{A.15}
\|\xi\|_{W^{1,2}(D^2)}+\|P\|_{W^{1,2}(D^2)}\le C(m)\ \|\Om\|_{L^2(D^2)}\quad ,
\ee
\item[iv)]
\be
\label{A.16}
\|\xi\|_{W^{2,2}(D^2)}+\|P\|_{W^{2,2}(D^2)}\le C(m)\ \|\Om\|_{W^{1,2}(D^2)}\quad .
\ee
\end{itemize}
\end{Lma}
{\bf Proof of lemma~\ref{lm-A3} $\Longrightarrow$ lemma~\ref{lm-A4}.}
Let $\Om$ in $L^2(D^2,so(m)\otimes\wedge^1{\R}^2)$ satisfying (\ref{A.8}). Introduce $\Om_k$
in $W^{1,2}(D^2,so(m)\otimes\wedge^1{\R}^2)$ converging strongly in $L^2$ to $\Om$.
Let $\xi_k$ and $P_k$ satisfying (\ref{A.13})...(\ref{A.16}) for $\Om_k$. Because of (\ref{A.15})
there exists a subsequence $\xi_{k'}$ and $P_{k'}$ converging weakly in $W^{1,2}$
to $\xi$ and $P$. Weak convergence in $W^{1,2}$ implies almost everywhere convergence
of $P_{k'}$ to $P$. Since $P_{k'}^tP_{k'}=id_{m}$, this equation passes to the limit and we have that
$P\in W^{1,2}(D^2,SO(m))$. Moreover, from Rellich compact embedding, $P_{k'}$ converges
strongly in every $L^q$ ($q<+\infty$) and $P_{k'}^t \Om_{k'} P_{k'}$ respectively 
$P_{k'}^t\nabla P_{k'}$ converge in distribution sense
to $P^t\Om P$ and respectively $P^t\nabla P$. Therefore (\ref{A.9}) is satisfied at the limit.
By continuity of the trace (\ref{A.10}) is also satisfied. Finally, by lower-semicontinuity of the
$W^{1,2}$ norm with respect to the weak $W^{1,2}$ convergence, we also obtain (\ref{A.11})
and lemma~\ref{lm-A3} is proved.

{\bf Proof of lemma~\ref{lm-A4}.}
We follow the strategy in \cite{Uhl}. We introduce the set
\be
\label{A.17}
{\mathcal U}_{\ep,C}=
\lf\{
\begin{array}{c}
\ds \Om\in W^{1,2}(D^2,so(m)\otimes\wedge^1{\R}^2)\mbox{ satisfying }\int_{D^2}|\Om|^2\le\ep\\[3mm]
 \mbox{ and for which there exists }\xi\in W^{2,2}(D^2,so(m))\\[3mm]
 \mbox{ and }P\in W^{2,2}(D^2,SO(m))\mbox{ solving }
 (\ref{A.13})\cdots(\ref{A.16})
\end{array}
\rg\}
\ee 
Following the previous argument we show that ${\mathcal U}_{\ep,C}$ is closed.
We now establish the following assertion:

{\bf Claim 1} {\it For any fixed $C$ there exists $\ep$ small enough, such that,
for any $\Om$ in ${\mathcal U}_{\ep,C}$
satisfying $\int_{D^2}|\Om|^2<\ep$ there exists a neighborhood of $\Om$ in $W^{1,2}$
included in ${\mathcal U}_{\ep,C}$.}

{\bf Proof of claim 1} Let $\Om\in {\mathcal U}_{\ep,C}$ satisfying $\int_{D^2}|\Om|^2<\ep$.
Let $\xi$ and $P$ satisfying (\ref{A.13})...(\ref{A.16}) for $\Om$. Following the arguments
in \cite{Uhl} (lemma 2.7 and 2.8), for every $\al>0$ we can find $\delta>0$ such that, for every $\la\in W^{1,2}(D^2,so(m)\otimes\wedge^1{\R})$ satisfying $\|\la\|_{W^{1,2}}\le\delta$, there exists $\xi_\la\in W^{2,2}(D^2,so(m))$ and 
$Q_\la\in W^{2,2}(D^2,SO(m))$ satisfying 
\be
\label{A.18}
\lf\{
\begin{array}{l}
\ds\nabla^\perp\xi_\la=Q_\la^{-1}\nabla Q_\la+Q^{-1}_\la(\nabla^\perp\xi+\la)Q_\la\quad\mbox{ in }
D^2\\[5mm]
\ds \xi_\la=0\quad\quad\mbox{ on }\p D^2
\end{array}
\rg.
\ee
and
\be
\label{A.19}
\|Q_\la-Id_m\|_{W^{2,2}}+\|\xi_\la-\xi\|_{W^{2,2}}\le\al\quad.
\ee
From (\ref{A.13}) and (\ref{A.18}) we then have
\be
\label{A.20}
\nabla^\perp\xi_\la=(PQ_\la)^{-1}\nabla(PQ_\la)+(PQ_\la)^{-1}(\Om+P\la P^{-1})PQ_\la\quad .
\ee
Since $P\in W^{2,2}$ the map $\la\rightarrow P\la P^{-1}$ and it's inverse $\zeta\rightarrow
P^{-1}\zeta P$ are continuous from $W^{1,2}$ into $W^{1,2}$ (using the fact that $W^{2,2}$ embedds 
in $W^{1,4}$ in 2 dimension and that $P\in SO(m)$). Therefore, for every $\beta>0$, there exists
$\eta>0$ such that for every $\zeta\in W^{1,2}(D^2,so(m)\otimes\wedge^1{\R})$ satisfying $\|\zeta\|_{W^{1,2}}\le\eta$,
there exists $R\in W^{2,2}(D^2,SO(m))$ and $\nu\in W^{2,2}(D^2,M_m({\R}))$ such that
\be
\label{A.21}
\lf\{
\begin{array}{l}
\ds\nabla^\perp\nu=R^{-1}\nabla R+R^{-1}(\Om+\zeta) R\quad\mbox{ in }D^2\\[5mm]
\ds \nu=0\quad\quad\mbox{ on }\p D^2
\end{array}
\rg.
\ee
Moreover we have
\be
\label{A.22}
\|R-P\|_{W^{2,2}}+\|\nu-\xi\|_{W^{2,2}}<\beta\quad.
\ee
Considering now $\beta<\ep^\frac{1}{2}$, since $\|P\|_{W^{1,2}}+\|\xi\|_{W^{1,2}}\le C\ep^\frac{1}{2}$
((\ref{A.15}) is satisfied for $\Om\in  {\mathcal U}_{\ep,C}$), we have
\be
\label{A.23}
\|R\|_{W^{1,2}}+\|\nu\|_{W^{1,2}}\le (C+1)\ep^\frac{1}{2}
\ee
In order to finish the proof of claim 1 it remains to establish (\ref{A.15}) and (\ref{A.16}), providing
that $\ep$ has been chosen small enough. This will be a consequence of the following lemma.
\begin{Lma}
\label{lm-A5}
There exist $C(m)>0$ $\ep$ and $\delta>0$ such that for every $P\in W^{2,2}(D^2,SO(m))$ and
 $\xi\in W^{2,2}(D^2,so(m))$ satisfying (\ref{A.13}) and (\ref{A.14}) for some $\Om\in
 W^{1,2}(D^2,so(m))$ satisfying $\int_{D^2}|\Om|^2\le\ep$, if
 \be
 \label{A.24}
 \|P\|_{W^{1,2}}+\|\xi\|_{W^{1,2}}\le\delta\quad,
 \ee
 then (\ref{A.15}) and (\ref{A.16}) are satisfied.
 \end{Lma}
 {\bf Proof of lemma~\ref{lm-A5}.}
 We first establish the critical estimate (\ref{A.15}). 
 
 (\ref{A.13}) and (\ref{A.14}) imply that
 $\xi$ solves the following elliptic PDE
 \be
 \label{A.25}
 \lf\{
 \begin{array}{l}
 \ds -\Delta\xi=\nabla P^t\cdot\nabla^\perp P+div(P^t \Om P)\quad\mbox{ in }D^2\quad ,\\[5mm]
 \ds \xi=0\quad\quad\mbox{ on }\p D^2
 \end{array}
 \rg.
 \ee
Using lemma~\ref{lm-A1} and standard elliptic PDE we have
\be
\label{A.26}
\|\nabla\xi\|_{L^2}\le C\ \|\nabla P^t\|_{L^2}\|\nabla P\|_{L^2}+ C\|\Om\|_{L^2}\quad .
\ee
Using the hypothesis that $\|\nabla P\|_{L^2}\le\delta$ we have that
\be
\label{A.27}
\|\nabla\xi\|_{L^2}\le C\ \delta\ \|\nabla P\|_{L^2}+ C\|\Om\|_{L^2}\quad .
\ee
From (\ref{A.13}) we have that
\be
\label{A.28}
\|\nabla P\|_{L^2}\le 2\|\nabla \xi\|_{L^2}+2\|\Om\|_{L^2}\quad.
\ee
Combining (\ref{A.27}) and (\ref{A.28}) we have then
\be
\label{A.29}
\|\nabla\xi\|_{W^{1,2}}\le 2C\ \delta\ \|\nabla \xi\|_{L^2}+ (C+2C\ \delta)\|\Om\|_{L^2}\quad .
\ee
Choosing then $2C\ \delta<1/2$ we obtain inequality (\ref{A.15}).

It remains to establish (\ref{A.16}). From (\ref{A.25}) again, using standard elliptic estimates
and the embedding of $W^{1,1}$ into $L^2$ in 2 dimensions,
we have
\be
\label{A.30}
\begin{array}{rl}
\ds\|\xi\|_{W^{2,2}}\le & C\|\nabla P^t\cdot\nabla^\perp P\|_{W^{1,1}}+C\|\Om\|_{W^{1,2}}\\[5mm]
 &\quad+
 C\|\nabla P^{t}\Om\|_{L^2}+C\|\Om\nabla P\|_{L^2}
 \end{array}
 \ee
Using Cauchy-Schwartz we have first
 \be
 \label{A.31}
 \| \nabla P^t\cdot\nabla^\perp P\|_{W^{1,1}}\le C\|P\|_{W^{2,2}}\|P\|_{W^{1,2}}\quad.
\ee 
Using the embedding of $W^{1,1}$ in $L^2$ and Cauchy-Schwartz we have
\be
\label{A.32}
\begin{array}{rl}
\|\nabla P^{t}\Om\|_{L^2}+\|\Om\nabla P\|_{L^2}\le &\ds C\|\nabla P^{t}\Om\|_{W^{1,1}}+\|\Om\nabla P\|_{W^{1,1}}\\[5mm]
 &\ds\le C\|P\|_{W^{2,2}}\|\Om\|_{L^2}+C\|\Om\|_{W^{1,2}}\|P\|_{W^{1,2}}\quad.
 \end{array}
\ee
Combining (\ref{A.30}), (\ref{A.31}) and (\ref{A.32}) we have then that
\be
\label{A.33}
\|\xi\|_{W^{2,2}}\le C(\delta +\ep^\frac{1}{2})\ \|P\|_{W^{2,2}} +C\|\Om\|_{W^{1,2}}\quad.
\ee
Using now (\ref{A.13}), we have
\be
\label{A.34}
\begin{array}{rl}
\ds\|\nabla P\|_{W^{1,2}}\le &\ds C\|\xi\|_{W^{2,2}}+C\||\nabla^\perp\xi| |\nabla^2 P|\|_{L^1}\\[5mm]
 &\ds+C\||\nabla^2\xi||\nabla P|\|_{L^1}
+C\|\Om\|_{W^{1,2}}+C\|\nabla P\Om\|_{W^{1,1}}\\[5mm]
 \le&\ds C(1+\delta)\|\xi\|_{W^{2,2}}+ C(1+\delta)\|\Om\|_{W^{1,2}}+C(\delta +\ep^\frac{1}{2})\ \|P\|_{W^{2,2}}
\end{array}
\ee
Combining (\ref{A.33}) and (\ref{A.34}), for $C(\delta+\ep^\frac{1}{2})<1/2$ we obtain
estimate (\ref{A.16}) and lemma~\ref{lm-A5} is proved.

\medskip

\noindent {\bf End of the proof of lemma~\ref{lm-A4}.}
Let $\Om$ in $W^{1,2}(D^2,so(m)\otimes\wedge^1{\R}^2)$ satisfying $\int_{D^2}|\Om|^2<\ep$ for $\ep$ for which
claim 1 holds.   We consider the path $\Om_t=\phi_t^\ast \Om$ where $\phi_t(x)=tx$ and $t\in[0,1]$. Since
 $\int_{D^2}|\Om_t|^2 =\int_{D^2_t}|\Om|^2$ is an increasing function of $t$ we have then
 a path among the elements in $W^{1,2}(D^2,so(m)\otimes{\R}^2)$ satisfying $\int_{D^2}|\Om_t|^2\le\ep$
connecting $0$ and $\Om$. Using the closedness of ${\mathcal U}_{\ep,C}$, the openess
property given by claim 1 and the fact that $0\in {\mathcal U}_{\ep,C}$, by the mean of a standard
continuity argument we obtain that $\Om$ is in ${\mathcal U}_{\ep,C}$ and lemma~\ref{lm-A4}
is proved.


\begin{thebibliography}{99}
\bibitem[Bet1]{Bet1} Bethuel, Fabrice "Un r\'esultat de r\'egularit\'e pour les solutions de l'\'equation de surfaces ˆ courbure moyenne prescrite". (French) [A regularity result for solutions to the equation of surfaces of prescribed mean curvature] C. R. Acad. Sci. Paris SŽr. I Math. 314 (1992), no. 13, 1003--1007.
\bibitem[Bet2]{Bet2} Bethuel, Fabrice ``Weak limits of Palais-Smale sequences for a class of critical functionals''.  Calc. Var. Partial Differential Equations  1  (1993),  no. 3, 267--310.
\bibitem[BeG1]{BeG1} Bethuel, F.; Ghidaglia, J.-M. Improved regularity of solutions to elliptic equations involving Jacobians and applications.  J. Math. Pures Appl. (9) 72 (1993), no. 5, 441--474.
\bibitem[BeG2]{BeG2} Bethuel, Fabrice; Ghidaglia, Jean-Michel "Some applications of the coarea formula to partial differential equations".  Geometry in partial differential equations, 1--17, World Sci. Publishing, River Edge, NJ, 1994.
\bibitem[BrC]{BrC} ÊBr\'ezis, Haim; Coron, Jean-Michel "Multiple solutions of $H$-systems and Rellich's conjecture".  Comm. Pure Appl. Math. 37 (1984), no. 2, 149--187.
\bibitem[CLMS]{CLMS} Coifman, R.; Lions, P.-L.; Meyer, Y.; Semmes, S. "Compensated compactness and Hardy spaces".  J. Math. Pures Appl. (9) 72 (1993), no. 3, 247--286.
\bibitem[Ev]{Ev} Evans Craig "Partial regularity for stationary harmonic maps into spheres" Arch. Rat. Mech.
Anal. 116 (1991), 101-113.
\bibitem[FMS]{FMS} Freire Alexandre, M\"uller Stefan and Struwe Michael "Weak compactness of wave maps
and harmonic maps" Ann. Inst. Henri Poincar\'e 15 (1998), no. 6, 725-754.
\bibitem[Gr1]{Gr1} Gr\"uter, Michael ``Conformally invariant variational integrals and the removability of isolated singularities.''
  Manuscripta Math.  47  (1984),  no. 1-3, 85--104.
\bibitem[Gr2]{Gr2}Gr\"uter, Michael ``Regularity of weak $H$-surfaces''.  J. Reine Angew. Math.  329  (1981), 1--15.
\bibitem[HiW]{HiW} Hildebrandt, Stefan; Widman, Kjell-Ove "Some regularity results for quasilinear elliptic systems of second order."  Math. Z. 142 (1975), 67--86.
\bibitem[Hei1]{Hei1}ÊHeinz, Erhard "Ein RegularitŠtssatz fŸr schwache Lšsungen nichtlinearer elliptischer Systeme". (German) Nachr. Akad. Wiss. Gšttingen Math.-Phys. Kl. II 1975, no. 1, 1--13. 
\bibitem[Hei2]{Hei2} Heinz, Erhard "†ber die RegularitŠt der Lšsungen nichtlinearer Wellengleichungen". (German) Nachr. Akad. Wiss. Gšttingen Math.-Phys. Kl. II 1975, no. 2, 15--26
\bibitem[Hei3]{Hei3} Heinz, Erhard "†ber die RegularitŠt schwacher Lšsungen nichtlinearer elliptischer Systeme". (German) [On the regularity of weak solutions of nonlinear elliptic systems] Nachr. Akad. Wiss. G\"ottingen Math.-Phys. Kl. II 1986, no. 1, 1--15. 
\bibitem[Hel]{Hel} ÊH\'elein, Fr\'ed\'eric "Harmonic maps, conservation laws and moving frames". 
 Cambridge Tracts in Mathematics, 150. Cambridge University Press, Cambridge, 2002.
 \bibitem[LaU]{LaU}Ladyzhenskaya, Olga A.; Ural'tseva, Nina N.  "Linear and quasilinear elliptic equations".
 Academic Press, New York-London 1968. 
\bibitem[LiR]{LiR} Lin, Fang-Hua; Rivi\`ere, Tristan "Energy quantization for harmonic maps".  Duke Math. J. 111 (2002), no. 1, 177--193.
\bibitem[MeR]{MeR} Meyer, Yves; Rivi\`ere, Tristan "A partial regularity result for a class of stationary Yang-Mills fields in high dimension."  Rev. Mat. Iberoamericana 19 (2003), no. 1, 195--219.
  \bibitem[Mo]{Mo}ÊMorrey, Charles B., Jr. "Multiple integrals in the calculus of variations". Die Grundlehren der mathematischen Wissenschaften, Band 130 Springer-Verlag New York, Inc., New York 1966
\bibitem[Mu]{Mu} M\"uller, Stefan  "Higher integrability of determinants and weak convergence in $L\sp 1$". 
J. Reine Angew. Math. 412 (1990), 20--34.
\bibitem[Sha]{Sha}Shatah, Jalal
Weak solutions and development of singularities of the ${\rm SU}(2)$ $\sigma$-model. 
Comm. Pure Appl. Math. 41 (1988), no. 4, 459--469.
 \bibitem[Ste]{Ste} Stein, Elias M.  "Singular integrals and differentiability properties of functions". Princeton Mathematical Series, No. 30 Princeton University Press, Princeton, N.J. 
 \bibitem[Ta1]{Ta1} Tartar, Luc "Remarks on oscillations and Stokes' equation.  Macroscopic modelling of turbulent flows" (Nice, 1984), 24--31, Lecture Notes in Phys., 230, Springer, Berlin, 1985
\bibitem[Uhl]{Uhl} Uhlenbeck, Karen K. ``Connections with $L\sp{p}$ bounds on curvature.''  Comm. Math. Phys.  83  (1982),
 no. 1, 31--42.
 \bibitem[We]{We} Wente, Henry C.  "An existence theorem for surfaces of constant mean curvature". 
J. Math. Anal. Appl. 26 1969 318--344.


 \end{thebibliography}
\end{document}